\documentclass[12pt,letterpaper,twoside]{amsart}

\usepackage[ngerman]{babel}
\usepackage{graphicx}
\usepackage{latexsym}
\usepackage{amsmath,amssymb}
\usepackage{ifpdf,stmaryrd}

\usepackage[DIV14,BCOR2cm]{typearea}

\theoremstyle{definition}
\newtheorem{Definition}{Definition}
\newtheorem{Theorem}[Definition]{Theorem}

\newtheorem{Lemma}[Definition]{Lemma}
\newtheorem{Example}[Definition]{Example}
\newtheorem{Vermutung}[Definition]{Conjecture}
\newenvironment{Beweis}{\noindent{Proof:}}{%
  \hspace*{\fill}$\Box$\par\vskip2ex}

\begin{document}


\title[On a conjecture concerning vertex-transitive graphs]{On a conjecture concerning vertex-transitive graphs}
\author{Tobias Ahsendorf}
\curraddr{Fakult\"at f\"ur Mathematik, Universit\"at Bielefeld, P.O.Box 100 131, D-33501 Bielefeld, Germany}
\email{tobias.ahsendorf@gmail.com}
\maketitle

\section*{Abstract}
In this article we define a minor relation, which is stronger than the classical one, but too strong to become a well-quasi-order on the class of finite graphs. Nevertheless, with this terminology we are able to introduce a conjecture, which would imply the Lovasz conjecture and give an interesting insight on the symmetry of vertex-transitive graphs, if true. Though it could become an approach to solve the Lovasz conjecture. These ideas were first introduced by the author in \cite{TA}.

\section{Introduction}
First of all, we have to introduce the classical minor relation with that we are able to establish the so-called homogeneous minor relation.
\begin{Definition}(see \cite{RDE}) \label{min}
Let $G$ and $H$ be graphs. We say $H$ is a \emph{minor} of $G$, in notation $H \preccurlyeq G$, when there is a partition $\{V_h\}_{h \in V(H)}$ of a subset of the vertex set of $G$ with $G[V_h]$ connected for each $h$, such that there is a $V_h$-$V_{h'}$ edge if $hh' \in E(H)$. 
\end{Definition}

There is a great theorem in minor theory, namely the Graph Minor Theorem, which was proved in twenty papers, published from 1983 till 2004 by Robertson and Seymour:
\begin{Theorem}(see \cite{RSXX}) \label{RS}The finite graphs are well-quasi-ordered by the minor relation $\preccurlyeq$, i.e. $\preccurlyeq$ is reflexive and transitive and for an infinite sequence of graphs $G_1,G_2,\dots$ there are $i,j$ with $i < j$ and $G_i \preccurlyeq G_j$.
\end{Theorem}

We actually do not need this Theorem here, it was only noted because this Theorem does not hold for the minor relation to be defined here.o

\section{The homogeneous minor relation}
In this section we want to introduce the already mentioned special minor relation compared to the relation of Definition \ref{min}, here the induced subgraphs partition sets of the vertex sets should be isomorphic. After this, we are going to establish some results concerning vertex-transitive graphs, pose the conjecture and prove that it would imply Lovasz conjecture.
\begin{Definition} Let $G,H,H'$ be graphs and $H'$ be connected. 
We call $(H,H')$ a \emph{homogeneous minor}, noted by $(H,H') \preccurlyeq_1 G$, if $H \preccurlyeq G$ and all partition sets $V_x$ can be chosen that way, that $G[V_x] \simeq H'$ holds.
\end{Definition}
If we write instead of $(H,H') \preccurlyeq_1 G$ just $H \preccurlyeq_1 G$ under the implicit assumption, that there exists a connected graph $H'$ with $(H,H') \preccurlyeq_1 G$, so we see, that we get a new minor relation. This minor relation is in opposite to the "'classical"' minor relation not a well-quasi-order on the class of the finite graphs. The following example is such a counterexample.
\enlargethispage{1 cm}
\begin{Example} Let $p_i$ denote $i$-th prime number and we regard the sequence $$C_3,\dots,C_{p_i},\dots$$ of cycle graphs of increasing and odd prime number order.  If $C_{p_i} \preccurlyeq_1 C_{p_j}$ for a pair $(i,j)$ with $1<i<j$, so the partition sets have to be paths of a certain length. Of course they are isomorphic, if the length of all these paths is equal. From this, we would obtain, that $p_i \,| \,p_j$, which is a contradiction to the fact, that $p_j$ is a prime number.
\end{Example}

Our next aim is to establish the already mentioned conjecture, that is connected to vertex-transitive graphs.

\begin{Definition}(cf. \cite{GR} 3.1)
\begin{itemize}
\item A graph $G$ is \emph{vertex-transitive}, if it has the property, that for $x,y \in V(G)$ there is an automorphism of $G$, i.e. an isomorphism from $G$ to $G$, such that this automorphism maps $x$ to $y$.
\item Let $G$ be a finite group und $C$ be a proper subset of $G$ with the properties, that the neutral element of $G$ is not contained in $C$ and that $x \in C$ implies $x^{-1} \in C$. Then the \emph{Cayley graph $X(G,C)$} is the graph consisting of the vertex set $G$ and the edge set $\{\,xy \, | \, xy^{-1} \in C \, \}$.
\end{itemize}
\end{Definition}
The Cayley graph has been defined in \cite{GR} for infinite groups, thus the graph itself would be infinite, but here we are just interested in the finite case. If we regard the second part of the preceding definition, we see among others, that $C$ was chosen that way, that no loops can appear in Cayley graphs, because the neutral element must not be an element of $C$, and that the Cayley graph can be regarded as undirected, because $C$ is closed under taking inverses. The most important property of Cayley graphs is its vertex-transitivity. We will prove it in our next Lemma. To do so we should establish the item of the automorphism group of a graph - this ist just the set of all automorphisms of the considered graph. This set is a (finite) group in the obvious way. 

\begin{Lemma}(cf. \cite{GR} Theorem 3.1.2) Let $G,C$ like in the preceding definition. Then the Cayley graph $X(G,C)$ vertex-transitive.
\end{Lemma}

\begin{Beweis} For every $g \in G$ we define a map $\rho_g: G \rightarrow G$ by $x \mapsto xg$. All of these maps are permutations of $G$. But they are automorphisms of $X(G,C)$ as well, because $(yg)(xg)^{-1} = ygg^{-1}x^{-1}=yx^{-1}$. Hence, $xg$ and $yg$ are adjacent, if $x$ and $y$ are adjacent. $\{\, \rho_g \,| \,g \in  G \, \}$ clearly build a subgroup of the automorphism group of $X(G,C)$, which is canonical isomorph to $G$. Now take two vertices of $X(G,C)$, say $g$ and $h$, then the automorphism $\rho_{g^{-1}h}$  maps $g$ to $h$.
\end{Beweis}

It is necessary to remark, that the converse does not hold, i.e. there are vertex-transitive graphs, which are not Cayleygraphs. For example, the Petersen graph is not a Cayley graph (siehe \cite{GR} Lemma 3.1.3). Surprisingly, the converse is true, if we just regard prime number orders.
\begin{Theorem}(cf. \cite{DM}) \label{DM1} Let $p$ be a prime number und $X$ a vertex-transitive graph with $p$ vertices. Then $X$ is a Cayley graph.
\end{Theorem}

We will not prove this Theorem here. In the just cited article \cite{DM} Maru\v{s}i\v{c} even proved, that this assertion is true even for $p^k$, if $k \leq 3$. Now we want to prove a result, which shows among others, that there is always a hamiltonian path in a connected vertex-transitive graphs, if the amount of vetices is a prime number. This has already been done in a similar form in \cite{DM}. 

\begin{Theorem} \label{DM2} Let $X$ be a vertex-transitive connected graph with $p$ vertices, where $p$ is a prime number. Then $X$ contains a hamiltonian path. Furthermore, if $p>2$, then $X$ contains a hamitonian cycle.
\end{Theorem}

\begin{Beweis} By Theorem \ref{DM1} there is a group $G$ with $|G|=p$ and a subset $C$ with $C^{-1} \subseteq C$ an the property, that the neutral element is not contained in $C$, such that $X=X(G,C)$ holds (or $X \simeq X(G,C)$, to be more precise). First of all, we regard $G$. By elementary group theory, we obtain that $G$ is isomorphic to $\mathbb{Z} / p \mathbb{Z}$. Hence, we can restrict to the last mentioned group. Because $C \neq \varnothing$, there is a $\overline{0} \neq k \in C$. Let us start at the vertex $\overline{0}$ and sum $k$ to the actual element, so we obtain after $p-1$ further additions of $k$ once again the residue class $\overline{0}$ and in the meantime we we visit each vertex exactly one time, if $p>3$ we obtain a hamiltonian cycle, in particular a hamiltonian path. But if $p=2$, we just have to consider $K_2$, and so we only have a hamiltonian path.
\end{Beweis}

Now we are going to formulate the wellknown Lovasz conjecture:
\begin{Vermutung}[Lovasz](cf. \cite{GR} 3. Notes) Every connected vertex-transitive graph contains a hamiltonian path.
\end{Vermutung}

Here we cannot exchange the expression of the hamiltonian path by the one of the hamiltonian cycle. This is easy to see by considering the Petersen graph.
From the soon established conjecture follows, together with the interplay of some already introduced results, the Lovasz conjecture.

\begin{Definition} Let $G$ be a vertex-transitive graph. Then, $G$ is in particular regular, with the \emph{valency} of such a graph we mean the valency of an arbitrary vertex of $G$.
\end{Definition}

\begin{Vermutung} \label{MV} Let $X$ be a vertex-transitive graph of valency $d$. Let $m,n \in \mathbb{N}$ with $mn=|X|$. ´
\begin{itemize}
\item If $(C_m,C_n)\preccurlyeq_1 X$, then $X$ contains a hamiltonian path. 
\end{itemize}
Now we ask additionally, that $|X|$ ist not a prime number and neither $0$ nor $1$.
\begin{itemize}
\item Then there is $m,n>1$ with $mn=|X|$, such that:
\begin{itemize}
\item $(C_m,C_n) \preccurlyeq_1 X$ or
\item There are two graph $X'$ and $X''$ with the following properties: $(X',X'') \preccurlyeq_1 X$ and $|X'|=m,|X''|=n$, such that $X'$ is vertex-transitive and connected, if $X$ is connected, and $\overline{K_m}$, if $X$ is not connected. \enlargethispage{1 cm}
Between two different vertices $x,y$ of $X''$, who have in $X''$ a degree $<d$, there is a hamiltonian path from $x$ to $y$. Additionally, if one regards the partition sets $V_x$ for all $x \in X'$ of $X'$, if $x_0 \in X'$ and $x_1,x_2 \in N_{X'}(x_0)$ with $x_1 \neq x_2$, so there exists for $v \in V_{x_0} \cap N_X(V_{x_1})$ a vertex $w \in V_{x_0}$ with $v\neq w$, such that there is a $\{w\}$-$V_{x_2}$-edge in $X$.
\end{itemize}
\end{itemize}
\end{Vermutung}
The principal part of this conjecture can be found in the second part. The first part is just established to conclude the Lovasz conjecture from this one, because, roughly spoken, the vertex-transitive graphs split in two classes because of the last part of the conjecture. If we obtain the first case, i.e.  $(C_m,C_n) \preccurlyeq_1 X$ for the regarded vertex-transitive graph $X$ with $m,n>1$ and $mn=|X|$, so we do not know, if there exists a hamiltonian path, therefore we need the first part of the conjecture. If we obtain the other case, we could, as later in the proof of the implication, work inductive, so we can directly prove the assertion without further assumptions.  If we would have proved the second part of the conjecture, then the only connected vertex-transitive graphs for which we have to verify the Lovasz conjecture, are those, which have a number of vertices, which is not $0$, $1$ or a prime number, and which does not admit a decomposition according to the last case of the second part of the conjecture \ref{MV}. Therefore this conjectrue could become an access to prove the Lovasz conjecture.\\
If it would reveal, that conjecture \ref{MV} is wrong, then it would be possible to weaken it in cooperation with the proof given below, such that this new conjecture would imply the Lovasz conjecture as well. Actually, a proof of the lower part of the just mentioned conjecture together with the Theorems \ref{DM1} and \ref{DM2} would give rise to an interesting classification of the class of vertex-transitive graphs. It is remarkable that we would get an insight on the symmetry of vertex-transive graph, if the conjecture is true. We could formulate this conjecture just for connected vertex-transitive graphs, because an unconnected vertex-transitive Graph $X$ consists of some copies of a connected vertex-transitive graph $X''$. In this case the last statement of the second part of the conjecture would fit with $X' = \overline{K_m}$, where $m$ denotes the multiplicity of $X''$ and $X''$ is as just defined. We include this case in the conjecture to obtain greatest possible generality. In order to prove the implication of the Lovasz conjecture we want to cite another Lemma, but we will not prove it here, which gives an anwer to the question on how much vertices of a connected vertex-transitive graph have to be removed at least to obtain a unconnected graph.

\begin{Lemma}(cf. \cite{GR}: Special case Theorem 3.4.2) \label{TS} We have to remove at least $\frac{2}{3}(k+1)$ from a connected vertex-transitive graph with valency $k$ to obtain a graph who is not connected.
\end{Lemma}

Now we are going to prove, that conjecture \ref{MV} together with Theorem \ref{DM2} and Lemma \ref{TS} implies the Lovasz conjecture. If we speak within the proof of "'the"' conjecture, we mean conjecture \ref{MV}.

\begin{Theorem} Assume conjecture \ref{MV}. Then the Lovasz conjecture is true.
\end{Theorem}
\enlargethispage{1 cm}
\begin{Beweis}
With $X$ we denote the regarded vertex-transitive connected graph with valency $d$. If $X$ has prime number order, the result follows by Theorem \ref{DM2}.
If $(C_m,C_n) \preccurlyeq_1 X$ for certain $m,n \in \mathbb{N}$ with $mn=|X|$, the claim follows by the first part of the conjecture. Especially, the statement is true for the empty graph and $K_1$. From now on we work inductive. Assume that $|X|=k$ and for all vertex-transitive connected graphs with less than $k$ vertices the claim has been verified. \\
Because of the considerations at the beginning of the proof we may assume, that $k>1$ is not a prime number and that $(C_m,C_n) \not\preccurlyeq_1 X$ for all $m,n \in \mathbb{N}$ with $mn =k$. 
Because of the conjecture we then have $(X',X'') \preccurlyeq_1 X$ with $X'$ and $X''$ having the properties described in the last part of the conjecture and that $X'$ is connected. Because $m<k$, $X'$ contains a hamiltionian path by induction. Now let $\{\, V_x \, | \, x \in V(X') \, \}$ be those partition sets of $X$ (or more precisely: of $V(X)$) with the property, that $X[V_x] \simeq X''$ for all $x \in V(X')$, such a choice is possible by assumption. Let $x_1\dots x_m$ be a hamiltonian path in $X'$. For all $1\leq i <m$ we can choose an edge $e_i$ in $X$ with the properties, that one vertex of this edge lies in $V_{x_i}$ and the other one is an element of $V_{x_{i+1}}$ and that $e_i \cap e_j = \varnothing$ for $i \neq j$. Again, this is possible because of the assumption. Let $v_1^{(1)}$ be a vertex in $V_{x_1}$ with $v_1^{(1)} \notin e_1$ and $d_{X[V_{x_1}]}(v_1^{(1)}) < d$. Analogical we define $v_m^{(2)}$ to be a vertex in $V_{x_m}$ with $v_m^{(2)} \notin e_{m-1}$ and $d_{X[V_{x_2}]}(v_m^{(2)})<d$. Both must exist, because otherwhise there would be just one vertex in $V_{x_1}$ (resp. $V_{x_m}$) with valency $<d$ in $X[V_{x_1}]$ (resp. $X[V_{x_m}]$). Thus we would obtain an unconnected graph by deleting the vertex in $V_{x_1} \cap e_1$ (resp. $V_{x_m} \cap e_{m-1}$). But since the valency of $X$ is greater or equal to $2$, this would contradict Lemma \ref{TS}. In the case, that $m>2$, the assertion is be easier to prove, otherwise $X[V_{x_1}] \simeq X[V_{x_m}] \simeq X[V_{x_2}] \simeq X''$ would not hold. Only in the case, that $m=2$, the statement is less obvious. \\
Now we are able to construct a hamiltonian path of $X$, who starts at $v_1^{(1)}$ and ends at $v_m^{(2)}$. Therefore let $e_i$ consist of $v_i^{(2)}$ and $v_{i+1}^{(1)}$ for all $1 \leq i <m$. Furthermore, it should hold, that $v_i^{(j)} \in V_{x_i}$ for $j=1,2$ and all $1 \leq i \leq m$. Because of the conjecture, there is for all $1 \leq i \leq m$ a hamiltonian path $P_i$ between $v_i^{(1)}$ and $v_i^{(2)}$ in $X[V_{x_i}]$. Now we can connect them, i.e. we build $P=P_1\dots P_m$, because $v_i^{(2)}$ and $v_{i+1}^{(1)}$  are adjacent by $e_i$ for all $1 \leq i <m$. Because every vertex of $X$ is contained in exactly one partition set, $P$ must be a hamiltonian path of $X$.
\end{Beweis}

\renewcommand{\refname}{References}

\vspace{0,3 cm}
\parindent=0pt
\end{document}